\documentclass[11pt]{article}
\usepackage{amsmath}
\usepackage{graphic x}

\topmargin\headheight
\usepackage{indentfirst}
\usepackage{amsfonts, amssymb,amsmath}
\newtheorem{t1}{Theorem}
\newtheorem{l1}{Lemma}
\newtheorem{proof}{Proof}

\begin{document}
	\baselineskip18pt
	
	\begin{center}\bf{\Large Fractional ordered Euler Riesz difference sequence spaces}
	\end{center}
	\begin{center}
		\centerline{\bf $Diptimayee Jena^{1}$, $Salila Dutta^{2,*}$ } 
	\centerline{ Department of Mathematics, Utkal University, India}
		 \centerline{1 E-mail: jena.deeptimayee@gmail.com}
		\centerline{2 E-mail:saliladutta516@gmail.com.}
		* - corresponding author
	\end{center}
\textbf{Abstract :}
{\textbf{}\small The main objective of this article is to introduce Euler-Riesz difference sequence spaces of fractional order ${\tau} $ along with infinite matrices. Some topological properties of these spaces are considered here along with the Schauder basis,  ${\alpha} -,\beta -$ and $\gamma -$duals of the spaces.\\}
	\centerline{}
{\bf 2010 Subject Classification:} 40A05, 46A45, 46B45, 46A35. \\

	\section{Introduction}
	For all real number $ \tau $ ,the gamma function  $\Gamma\left({\tau} \right)$ is expressed as 
	 \begin{equation}
		\Gamma\left({\tau} \right)=\displaystyle\int_0^\infty e^{-t} \ t^{\tau-1}  \, dt
	\end{equation} , which is an improper integral and it satisfies the following properties: 
     
     1. $\Gamma\left({n+1} \right)=n!$,  $ n \in N$
     
     2. $\Gamma\left({n+1} \right)=n \Gamma\left({n} \right)$
     for each real number $n\not\in \{0,-1,-2,-3,....\}$.\
    
    \noindent It is obvious by previous knowledge based on sequence spaces the space of all real valued sequences are denoted by $\omega $. The spaces $c_{0} ,c $ and $l_{\infty }$  are of all  null, convergent and bounded sequences respectively.\
     
      Maddox \cite{madx1,madx2,madx3} introduced $c_{0}\left(p\right) ,c\left(p\right) $ and $l_{\infty }\left(p\right) $ for $p=\left(p_{k}\right) $ a bounded sequence of strictly positive real numbers, as:
     \[\begin{array}{l} {c_{0}\left(p\right) =\left\{\zeta=\left(\zeta_{k} \right)\in \omega :{\mathop{\lim }\limits_{k\to \infty }} \left| \zeta_{k} \right|^{p_{k}} =0 \right\}}, \\ 
     	{c\left(p\right) =\left\{\zeta=\left(\zeta_{k} \right)\in \omega :{\mathop{\lim }\limits_{k\to \infty }} \left| \zeta_{k} -l\right|^{p_{k}} =0 \ for some \ l\in R \right\}}, \\ 
     	{l_{\infty }\left(p\right) =\left\{\zeta=\left(\zeta_{k} \right)\in \omega :{\mathop{\sup }\limits_{k\in N}} {\left| \zeta_{k} \right|}^{p_{k} } <\infty \right\}} \end{array}\]  \
    	and these are complete paranormed sequence spaces with paranorm  
    $ g\left(x\right) = {\mathop{\sup }\limits_{k\in N}}\ {\left| x_{k} \right|}^{{p_{k} }/M } $ and  where $ M=max \{1,\displaystyle \sup_{k} \ p_{k}\} $. The difference sequence spaces $X\left(\Delta \right)=\left\{x=\left(x_{k} \right):\Delta \left(x_{k} \right)\in X\right\}$ for $X=\left\{c_{0}, c, l_{\infty }  \right\}$  was introduced in 1981 by K1zmaz \cite{kizm15}. Then the difference sequence spaces attracted the attention of several authors \cite{basa6,colak1,mur1} in different directions.  For a proper fraction ${\tau}$, Baliarsingh \& Dutta  (\cite{bali4,bali5,bali6,bali7,dutt10,dutt11}) they introduced the difference operator $\Delta ^{\left( {\tau}\right) } $ as
    \begin{equation} \label{1.2} 
    	\Delta ^{\left( {\tau}\right)  } x_{k} =\sum _{i}\left(-1\right)^{i} \frac{\Gamma \left({{\tau} }+1\right)}{i!\Gamma ({{\tau} }-i+1)} x_{k-i}   
    \end{equation} 
    
    \begin{equation} \label{1.3} 
    	\Delta ^{-\left( {\tau}\right) } x_{k} =\sum _{i}\left(-1\right)^{i} \frac{\Gamma \left({{-\tau} }+1\right)}{i!\Gamma ({-{\tau} }-i+1)} x_{k-i}   
    \end{equation}   
    
    Throughout the paper it is assumed that the series of fractional difference operators are convergent. It is also appropriate to express the difference operator and its inverse as follows:\
    
    \[\left(\Delta ^{\left( {\tau}\right)  } \right)_{nk} =\left\{\begin{array}{l} {\left(-1\right)^{n-k} \frac{\Gamma \left({{\tau} }+1\right)}{\left(n-k\right)!\Gamma ({{\tau} }-n+k+1)} \qquad if\, 0\le k\le n} \\ {0 \qquad \qquad \qquad \qquad \qquad \qquad if\, k>n} \end{array}\right. \] \
    
    \[\left(\Delta ^{\left( {-\tau}\right)  } \right)_{nk} =\left\{\begin{array}{l} {\left(-1\right)^{n-k} \frac{\Gamma \left({-{\tau} }+1\right)}{\left(n-k\right)!\Gamma ({-{\tau} }-n+k+1)} \qquad if\, 0\le k\le n} \\ {0 \qquad \qquad \qquad \qquad \qquad \qquad if\, k>n} \end{array}\right. \] \
    Basar and Braha \cite{bara3} introduced Euler-Cesaro difference sequence spaces  $\check{c}$ , $\check{c_{0}}$, $\check{l_\infty}$ of null, convergent and bounded sequences respectively and then Ellidokuzoglu and Demiriz \cite{demri2} introduced Euler-Riesz difference sequence spaces  ${ \left[{ c_{0}} \right] _{e,r}}$ , ${ \left[ c \right] _{e,r}}$ , ${\left[ {l_\infty} \right] _{e,r}}  $. Baliarsingh and Dutta \cite{bali4} introduced the fractional difference operators. Now our interest is to introduce Euler-Riesz difference sequence spaces of fractional order.\
    
    \noindent Now we introduce the spaces $ {\left[ c_{0} , \Delta^{\left(\tau \right)} , p \right]}_{er} $ , $ {\left[ c , \Delta^{\left(\tau \right)} , p \right]}_{er} $ and $ {\left[ l_\infty , \Delta^{\left(\tau \right)} , p \right]}_{er} $ by using the product of the Euler mean $ E_{1} $ and Riesz mean $ R_{q} $ with fractional operator $ \Delta^{\left( \tau \right)} $ .We prove certain topological properties of these spaces along with $ \alpha- , \beta- , \gamma- $ duals.\
    
    \section{Main Results}
    
    \noindent Here we introduce the matrix  $ {\tilde{B} \left( \Delta^{\left( \tau \right)}\right)} ={\tilde{B}^{\left( \tau \right)}} $ by the product of Euler-Riesz matrix $\tilde{B}$ \cite{demri2} and fractional ordered difference operator  $ \Delta^{\left( \tau \right)}$ \cite{bali4},  and obtain its inverse, where\
     \begin{multline} \label{2.1} 
    \left({\tilde{B}^{\left( \tau\right)}} \right)_{nk} =\left\{\begin{array}{l} {\sum _{i=k}^{n}(-1)^{i-k} }  \left(\begin{array}{l} {n} \\ {i} \end{array}\right)\frac{\Gamma \left({\tau} +1\right)}{\left(i-k\right)!\Gamma \left({\tau} -i+k+1\right)} \ \frac{q_{i}}{ 2^{n}  Q_{n}} , { \qquad if\, 0\le k\le n } \\ {0,\qquad \qquad \qquad \qquad \qquad \qquad\qquad \qquad \qquad \qquad if\, k>n} \end{array}\right. 
    \end{multline} 

    Equivalently we may write\
    \begin{multline}\label{*}
      \left({\tilde{B}^{\left( \tau\right)}} \right) =\
    	\left(%
    	\begin{array}{ccccccc}
    		{\frac{1}{2} } & 0 & 0 & 0 & \dots \\
    		({\frac{2 q_{1}- \tau q_{2}}{2^2 Q_{2} } })  &  ({\frac{q_{2} }{2^2 Q_{2}} })   & 0 & 0 & \dots  \\
    		({\frac{3q_{1}-3{\tau} q_{2} +\frac{{\tau}({\tau}-1)}{2!} q_{3} }{2^3 Q_{3} } })  &  ({\frac{3 q_{2}- \tau q_{3}}{2^3 Q_{3} } })  &  ({\frac{q_{3}}{2^3 Q_{3}}})  & 0&\dots  \\
    		\vdots  & \vdots  & \vdots  &\vdots &\ddots
    	\end{array}%
    	\right)
    \end{multline}\\

     By simple calculation the inverse of the matrices $\Delta^ {\left(\tau \right)} , \tilde{B}$ and $\tilde{B}^{\left( \tau\right)}$  can be obtained as given in the following lemma. 
   
    \begin{l1} \label{1.1}\cite{bali7}
    	\noindent The inverse of fractional difference operator $\Delta ^{\left( {\tau}\right) } $ is given as
    	\[\left(\Delta ^{-\left({\tau}\right)} \right)_{nk} =\left\{\begin{array}{l} {\left(-1\right)^{n-k} \frac{\Gamma \left(-{\tau} +1\right)}{\left(n-k\right)!\Gamma (-{\tau} -n+k+1)} \qquad if\, 0\le k\le n} \\ {0\qquad \qquad \qquad \qquad  \qquad \qquad if\, k>n} \end{array}\right. \] 
    	
    \end{l1}
\begin{l1}\label{2}\cite{demri2}
  	\noindent The inverse of the Euler-Riesz matrix $\tilde{B}$ is given as
  	\[\left(\tilde{B} \right)_{nk}^{-1} =\left\{\begin{array}{l} {(-1)^{n-k} \ \  \left(\begin{array}{l} {n} \\ {k} \end{array}\right) {\frac{Q_{k}{2^k}}{q_{n}}}, \qquad  \qquad if\, 0\le k\le n} \\ {0,  \qquad    \qquad \qquad \qquad \qquad \qquad if\, k>n} \end{array}\right. \] 
  \end{l1}
\begin{t1}\label{4}
	\noindent The inverse of the fractional ordered Euler-Riesz matrix  
	$  {\left({\tilde{B}^{\left(\tau \right)}}\right)}_{nk} $ is written as $ {\left({\tilde{B}^{\left(\tau \right)}}\right)}_{nk}^{-1} $ and given by	
	\[{\left({\tilde{B}^{\left(\tau \right)}}\right)}_{nk}^{-1} =\\
	\left\{\begin{array}{l} { \sum _{j=k}^{n}\left(-1\right)^{n-k}  \left(\begin{array}{l} {j} \\  {k} \end{array}\right) \frac{\Gamma \left({-\tau} +1\right)}{\left(n-j\right)!\Gamma \left({-\tau}-n+j+1\right)} \frac{2^k Q_{k}}{q_{j}} , \qquad \qquad   if\, 0\le k\le n } \\ {0, \qquad \qquad \qquad \qquad \qquad \qquad \qquad \qquad  \qquad  \qquad  \qquad  if\, k>n} \end{array}\right. \] 	
\end{t1}

\begin{proof}
This theorem can be proved using lemma \ref{1.1}, lemma \ref{2}  and hence omitted.
\end{proof}

 \noindent For a positive real number $\tau $ ,  we now introduce the classes of fractional ordered Euler-Riesz difference sequence spaces  $ {\left[ c_{0} , \Delta^{\left(\tau \right)} , p \right]}_{er} $ , $ {\left[ c , \Delta^{\left(\tau \right)} , p \right]}_{er} $ and $ {\left[ l_\infty , \Delta^{\left(\tau \right)} , p \right]}_{er} $ by
 
 \begin{tiny}
 	\[\begin{array}{l} {\left[ c_{0} , \Delta^{\left(\tau \right)} , p \right]}_{er} 
 {=\left\{x=\left(x_{k} \right)\in w:{\mathop{\lim }\limits_{m\to \infty }} \left|  \sum _{j=0}^{m}\sum _{i=j}^{m} \left(-1\right)^{i-j}  \left(\begin{array}{l} {m} \\ {i} \end{array}\right)
 	\frac{\Gamma \left({\tau} +1\right)}{\left(i-j\right)!\Gamma \left({\tau} -i+j+1\right)}\frac{q_{i} x_{j}}{ 2^{m}  Q_{m}}\right| ^{p_{k} }
 	=0  \right\}} \end{array}\]
 \[\begin{array}{l}  {\left[ c , \Delta^{\left(\tau \right)} , p \right]}_{er}  {=\left\{x=\left(x_{k} \right)\in w:{\mathop{\lim }\limits_{m\to \infty }} \left|  \sum _{j=0}^{m}\sum _{i=j}^{m} \left(-1\right)^{i-j}  \left(\begin{array}{l} {m} \\ {i} \end{array}\right)
 		\frac{\Gamma \left({\tau} +1\right)}{\left(i-j\right)!\Gamma \left({\tau} -i+j+1\right)}\frac{q_{i} x_{j}}{ 2^{m} Q_{m}}\right| ^{p_{k} }
 		exists  \right\}} \end{array}\] 
 \[\begin{array}{l} {\left[ l_\infty , \Delta^{\left(\tau \right)} , p \right]}_{er} {=\left\{x=\left(x_{k} \right)\in w:{\mathop{\sup }\limits_{m}} \left|  \sum _{j=0}^{m}\sum _{i=j}^{m} \left(-1\right)^{i-j}  \left(\begin{array}{l} {m} \\
 			 {i} \end{array}\right)
 		\frac{\Gamma \left({\tau} +1\right)}{\left(i-j\right)!\Gamma \left({\tau} -i+j+1\right)}\frac{q_{i} x_{j}}{ 2^{m}  Q_{m}}\right| ^{p_{k} }
 		<\infty \right\}}    \end{array}\] 
 
\end{tiny}
Above spaces can be written as :\\
$ {\left[ c_{0} , \Delta^{\left(\tau \right)} , p \right]}_{er} 
	=\left(c_{0}\left( p \right)  \right)_{\left(\tilde{B}^{\left( \tau\right) }  \right)}  , \left[ c ,\ \Delta^{\left(\tau \right)} , p \right]_{er} 
{=\left(c\left( p \right)  \right)_{ \left( \tilde{B}^{\left( \tau\right) } \right) }, {\left[ l_\infty ,\ \Delta^{\left(\tau \right)} ,p \right]}_{er}=\left(l_\infty\left( p \right) \right)_{\left(\tilde{B}^{\left( \tau\right) }\right)} }.$
The above sequence spaces generalize many of the known sequence spaces as particular cases which are as follows:\\
${\left( i\right)} $. For ${ \tau ={0} }$ and ${ p=\left(  p_{k}\right) = e}$ , the above classes reduce to ${ \left[{ c_{0}} \right] _{e,r}}$ , ${ \left[ c \right] _{e,r}}$ , ${\left[ {l_\infty} \right] _{e,r}}  $ introduced by Ellidokuzoglu and Demiriz.\cite{demri2}\\
${\left( ii\right)} $. For ${ \tau ={0} }$ and ${ p=\left( p_{k}\right) = e} , { q=\left( q_{k}\right) = e}$ , the above classes reduce to the sequence spaces $\check{c}$ , $\check{c_{0}}$, $\check{l_\infty}$ studied by Basar and Braha \cite{bara3}.\\
\noindent  Now with $\tilde{B}^{\left({\tau}\right)}$ - transform of $x={\left(x_{k}\right)}$ we define the sequence $y=\left(y_{k} \right)$ as follows :
\begin{equation} \label{3.1} 
	\begin{array}{l} {y_{k} =\left( \tilde{B}^{\left({\tau}\right)} x \right)_{k} } {=\sum _{j=0}^{k}\sum _{i=j}^{k}\left(-1\right)^{i-j}  \left(\begin{array}{l} {k} \\ {i} \end{array}\right)\frac{\Gamma \left({\tau} +1\right)}{\left(i-j\right)!\Gamma \left({\tau} -i+j+1\right)} \frac{q_{i} x_{j}}{ 2^{k}  Q_{k}} } \end{array} ,
\end{equation} 
By a straightforward calculation $\left( \ref{3.1}\right)$ it can be obtained that 
\begin{equation} \label{3.2} 
	\begin{array}{l} {x_{k} =\left(\tilde{B}^{\left({-\tau}\right)} y \right)_{k} } {=\sum _{i=0}^{k}\sum _{j=i}^{k}\left(-1\right)^{k-i}  \left(\begin{array}{l} {j} \\ {i} \end{array}\right) \frac{\Gamma \left({-\tau} +1\right)}{\left(k-j\right)!\Gamma \left({-\tau}-k+j+1\right)} \frac{2^i Q_{i} y_i}{q_{j}} } \end{array} ,
\end{equation} 
\begin{l1}\label{3.5}
	\noindent The operator $\tilde{B}^{\left({\tau}\right)}:  w \ {\rightarrow} \ w $ is linear.
\end{l1}
\begin{proof}
	The proof is omitted as it can be easily obtained. 
\end{proof}
Remark: $\tilde{B}^{\left({\tau}\right)} . \tilde{B}^{\left({-\tau}\right)} \cong \tilde{B}^{\left({-\tau}\right)} .\tilde{B}^{\left({\tau}\right)} \cong I $ , where I is an identity matrix.

\section{Topological structure} 
\noindent This section deals with some interesting topological results of the spaces $ {\left[ c_{0} , \Delta^{\left(\tau \right)} , p \right]}_{er} $, $ {\left[ c , \Delta^{\left(\tau \right)} , p \right]}_{er} $ and $ {\left[ l_\infty , \Delta^{\left(\tau \right)} , p \right]}_{er} $ .
\begin{t1}\label{th3.1}
 	\noindent $ {\left[ c_{0} , \Delta^{\left(\tau \right)} , p \right]}_{er} $ is a paranormed space with the  paranorm \[ {g}_{\tilde{B}^{\left( \tau\right)}} \left( x\right){= {\mathop{\sup }\limits_{k\in N}}\left|{ \left(\left(  \tilde{B}^{\left({\tau}\right)}\right)x\right)}_k\right| ^{\frac{p_{k}}{M}}}\]
 	\begin{equation}\label{pnorm1}
 	={\mathop{\sup }\limits_{n}} \left|  \sum _{j=0}^{k}\sum _{i=j}^{k} \left(-1\right)^{i-j}  \left(\begin{array}{l} {k} \\
 			{i} \end{array}\right)
 		\frac{\Gamma \left({\tau} +1\right)}{\left(i-j\right)!\Gamma \left({\tau} -i+j+1\right)}\frac{q_{i}}{ 2^{k}  Q_{k}} {x_{j}}\right| ^{\frac{p_k}{M}},
 			\end{equation}
 		 where $ 0\ <{p_k} \leq H\ < \infty , H={\mathop{\sup }\limits_{k}}\ {p_k} , M = max(1,H) , h={\mathop{\inf }\limits_{k}}\ {p_k}.$
 \end{t1}
\begin{proof}
\noindent The theorem will be proved only for $ {\left[ c_{0} , \Delta^{\left(\tau \right)} , p \right]}_{er} $.\\
 	 Assumed that $h > 0$ , then $ {g}_{\tilde{B}^{\left( \tau\right) }} \left( \theta \right)=0 $ and $ {g}_{\tilde{B}^{\left( \tau\right) }} \left( -x\right)={g}_{\tilde{B}^{\left( \tau\right) }} \left( x\right) $. To prove the linearity of $ {g}_{\tilde{B}^{\left( \tau\right) }} \left( x\right) $, we consider two sequences $x= {\left( x_k\right)} , y= {\left( y_k\right)} \in {\left[ c_{0} , \Delta^{\left(\tau \right)} , p \right]}_{er}$ and any two scalars ${{\beta_1} , {\beta_2} \in \mathbb{R} }$ .
 	 Since $\tilde{B}^{\left({\tau}\right)}$ is linear we get ,\\
 	 $ {g}_{\tilde{B}^{\left({\tau}\right)}} \left( {{\beta_1}x + {\beta_2}y} \right) ={\mathop{\sup }\limits_{n}} \left|  \sum _{j=0}^{k}\sum _{i=j}^{k} \left(-1\right)^{i-j}  \left(\begin{array}{l} {k} \\
 	 	{i} \end{array}\right)
 	 \frac{\Gamma \left({\tau} +1\right)}{\left(i-j\right)!\Gamma \left({\tau} -i+j+1\right)}\frac{q_{i}}{ 2^{k}  Q_{k}} {\left( {{\beta_1}{x_j} + {\beta_2}{y_j}}\right) }\right| ^{\frac{p_k}{M}} \\
 	  \leq {max\{ 1,\left| \beta_1\right| \} {\mathop{\sup }\limits_{k\in N}}\left|{ \left(\left(  \tilde{B}^{\left({\tau}\right)}\right)x\right)}_k\right| ^{\frac{p_{k}}{M}} + max\{ 1,\left| \beta_2\right| \} {\mathop{\sup }\limits_{k\in N}}\left|{ \left(\left(  \tilde{B}^{\left({\tau}\right)}\right)y\right)}_k\right| ^{\frac{p_{k}}{M}}} \\
 	  {= max\{ 1,\left| \beta_1\right| \} \ {{g}_{\tilde{B}^{\left({\tau}\right)}} \left( x\right)} + max\{ 1,\left| \beta_2\right| \} \ {{g}_{\tilde{B}^{\left({\tau}\right)}} \left( y\right)} }  $\\
 	  Hence ${g}_{\tilde{B}^{\left({\tau}\right)}}$ is subadditive. i.e. ${{g}_{\tilde{B}^{\left({\tau}\right)}} \left( {x+y}\right)} \leq {{g}_{\tilde{B}^{\left({\tau}\right)}} \left( x\right)} + {{g}_{\tilde{B}^{\left({\tau}\right)}} \left( y\right)}$ , for all $x,y \in {\left[ c_{0} , \Delta^{\left(\tau \right)} , p \right]}_{er} $. Now consider $ \{ u^n \} $ is a sequence of points in ${\left[ c_{0} , \Delta^{\left(\tau \right)} , p \right]}_{er}$ then $ {{g}_{\tilde{B}^{\left({\tau}\right)}} \left( {u^n} - u\right)} \rightarrow {0}$ and $\left( {\lambda_n} \right) $ is a sequence of scalars such that $ \lambda_n \rightarrow \lambda $ as $ {n} \rightarrow \infty $. By using the subadditivity of ${g}_{\tilde{B}^{\left({\tau}\right)}}$ , we get
 	     	$ {g}_{\tilde{B}^{\left({\tau}\right)}} \left( {u^n} \right) \leq {g}_{\tilde{B}^{\left({\tau}\right)}} \left( u\right)+{g}_{\tilde{B}^{\left({\tau}\right)}} \left( {u^n} - u\right) $.\\
 	  \noindent Since $ \{{g}_{\tilde{B}^{\left({\tau}\right)}} \left( {u^n}\right) \} $ is bounded , we have \\
 	   $ {g}_{\tilde{B}^{\left({\tau}\right)}} \left( { \lambda_n u^n} - \lambda u\right)  ={\mathop{\sup }\limits_{k}} \left|  \sum _{j=0}^{k} \left[ \sum _{i=j}^{k} \left(-1\right)^{i-j}  \left(\begin{array}{l} {k} \\
 	   	{i} \end{array}\right)
 	   \frac{\Gamma \left({\tau} +1\right)}{\left(i-j\right)!\Gamma \left({\tau} -i+j+1\right)}\frac{q_{i}}{ 2^{k}  Q_{k}}\right]  {\left( {{\lambda_n}{u_j}^n + {\lambda}{u_j}}\right) }\right| ^{\frac{p_k}{M}} \\
 	   \leq  {\left| {\lambda_n - \lambda}\right|}^{\frac{p_k}{M}} \  {g}_{\tilde{B}^{\left({\tau}\right)}} \left( {u^n} \right) + {\left| {\lambda}\right|}^{\frac{p_k}{M}} \  {g}_{\tilde{B}^{\left({\tau}\right)}} \left( {u^n} - u \right) \rightarrow 0 \ as\  {n \rightarrow \infty} $.\\
 	   Hence it shows that the scalar multiplication of ${g}_{\tilde{B}^{\left({\tau}\right)}} \left( {x}\right)$ is continuous and ${g}_{\tilde{B}^{\left({\tau}\right)}} \left( {x}\right)$ is a paranorm of the space ${\left[ c_{0} , \Delta^{\left(\tau \right)} , p \right]}_{er}$ and proof for other spaces can be done using similar techniques.  
 \end{proof} 
  \begin{t1}\label{th3.2}
  	\noindent The sequence space ${\left[ c_{0} , \Delta^{\left(\tau \right)} , p \right]}_{er}$ is a complete linear space paranormed by ${g}_{\tilde{B}^{\left({\tau}\right)}} \left( {x} \right)$ .
  	\end{t1}
  \begin{proof}
  	The proof is a routine verification and hence omitted. 	
  \end{proof}
 \begin{t1}\label{th3.3}
	\noindent $ {\left[ c_{0} , \Delta^{\left(\tau \right)} , p \right]}_{er} $, $ {\left[ c , \Delta^{\left(\tau \right)} , p \right]}_{er} $ and $ {\left[ l_\infty , \Delta^{\left(\tau \right)} , p \right]}_{er} $ 
	are linearly isomorphic to $c_{0} {\left(p \right)}  ,\, c {\left(p \right)},\, l_{\infty } {\left(p \right)} $  where $0\ <{p_k} \leq H\ < \infty$,\  respectively.
	\begin{proof}

	\noindent Now define a mapping $F:{\left[ l_\infty , \Delta^{\left(\tau \right)} , p \right]}_{er} \to l_{\infty} {\left(p \right)}  $ by      $x\to y=Fx$.
	\noindent Clearly , $F$ is a linear transformation. If $Fx=\theta \ then \ x=\theta $, so  $F$ is one-one. \
		
		\noindent Let $y\in l_{\infty}{\left(p \right)} $ , define a sequence $x=\left(x_{k} \right)$ in $\left( \ref{3.1}\right)$  as
			\[x_{k} =\sum _{i=0}^{k} \sum _{j=i}^{k}\left(-1\right)^{k-i}  \left(\begin{array}{l} {j} \\ {i} \end{array}\right)\frac{\Gamma \left(-{\tau} +1\right)}{\left(k-j\right)!\Gamma \left(-{\tau} -k+j+1\right)}  \frac{2^i Q_{i}}{q_{j}} y_{i}  \] 
		Then \[ {g}_{\tilde{B}^{\left({\tau}\right)}} \left( x\right) 
		={\mathop{\sup }\limits_{n}} \left|  \sum _{j=0}^{k}\sum _{i=j}^{k} \left(-1\right)^{i-j}  \left(\begin{array}{l} {k} \\
			{i} \end{array}\right)
		\frac{\Gamma \left({\tau} +1\right)}{\left(i-j\right)!\Gamma \left({\tau} -i+j+1\right)}\frac{q_{i}}{ 2^{k}  Q_{k}} {x_{j}}\right| ^{\frac{p_k}{M}} \] 
  		\[ = \mathop{\sup }\limits_{k \in N} \left|\sum _{j=0}^{k} {\delta_{kj} \ y_{j}}\right|^{\frac{p_k}{M}} ={\mathop{\sup }\limits_{k \in N}} \left| y_{k} \right|^{\frac{p_k}{M}} \ < \infty \]
 		\[ where \  {\delta_{kj}} = \left\lbrace  \begin{array}{l} { 1 , \qquad if\, k = j } \\ { 0 , \qquad if\, k \neq j } \end{array}\right. \]  
	Thus, $ x \in {\left[ l_\infty , \Delta^{\left(\tau \right)} , p \right]}_{er}  $ and $F$ is a linear bijection and paranorm preserving. Hence the spaces ${\left[ l_\infty , \Delta^{\left(\tau \right)} , p \right]}_{er}$ and $l_\infty\left( p\right) $ are linearly isomorphic.
		 \noindent i.e.  ${\left[ l_\infty , \Delta^{\left(\tau \right)} , p \right]}_{er} \cong l_\infty\left( p\right) $. The proof for other spaces can be obtained in a similar manner.
	\end{proof}
 \end{t1}
    \section{Basis for the spaces}
    \noindent In this section the Schauder basis \cite{madx3} for $ {\left[ c_{0} , \Delta^{\left(\tau \right)} , p \right]}_{er} $, $ {\left[ c , \Delta^{\left(\tau \right)} , p \right]}_{er} $ are constructed. 
    
      \begin{t1}\label{th4.1}
      For $0\ <{p_k} \leq H\ < \infty$ , let $ \mu_{k} \left( q\right){= { \left(\left(  \tilde{B}^{\left({\tau}\right)}\right)x\right)}_k} $. For $k \in N_{0}$ define $ {b^{\left( k\right)}} {\left( q\right)} = {\left\lbrace {{b_{n}}^{\left( k\right)}} {\left( q\right)}\right\rbrace}_{n \in N_{0}} $ by
      \begin{multline}
      	{\left\lbrace {{b_{n}}^{\left( k\right)}} {\left( q\right)}\right\rbrace} =\left\{\begin{array}{l} {\sum _{j=k}^{n}(-1)^{n-k} }  \left(\begin{array}{l} {j} \\ {k} \end{array}\right)\frac{\Gamma \left({-\tau} +1\right)}{\left(n-j\right)!\Gamma \left({-\tau} -n+j+1\right)} \ \frac{2^{k}  Q_{k}}{q_{j} } , { \qquad if\, 0\le k \le n } \\\ {\frac{2^{n} Q_{n}}{q_{n}}, \qquad \qquad \qquad \qquad \qquad \qquad\qquad\qquad \qquad if\, k=n} \\\ \   {0,\qquad \qquad \qquad \qquad \qquad \qquad\qquad \qquad \qquad \qquad if\, k>n} \end{array}\right. 
      \end{multline} 
  \noindent$(i)$  $\left\lbrace {{b_{n}}^{\left( k\right)}} {\left( q\right)}\right\rbrace$ is a basis for ${\left[ c_{0} , \Delta^{\left(\tau \right)} , p \right]}_{er} $ and each $x\in {\left[ c_{0} , \Delta^{\left(\tau \right)} , p \right]}_{er} $ and $ x $ has unique representation
  		\[x=\sum _{k}\mu _{k}{\left( q\right)} {{b_{n}}^{\left( k\right)}} {\left( q\right)} \] .
 \noindent$(ii)$  $\left\lbrace {\left(\tilde{B}^{\left( -\tau \right) } \right){e}}, {{b_{n}}^{\left( k\right)}} {\left( q\right)}\right\rbrace$ is a basis for ${\left[ c , \Delta^{\left(\tau \right)} , p \right]}_{er} $, and each $x\in {\left[ c , \Delta^{\left(\tau \right)} , p \right]}_{er} $ and $ x $ has unique  representation
 	\[x= {le} + \sum _{k} \left( {\mu_{k}-l}\right)   {b^{\left( k\right)}}, \  \ l=\lim\limits_{k \rightarrow \infty} {\mu _{k}} \ . \]   	
  \end{t1}
   \begin{proof}
   	\noindent $(i)$  By the definition of ${\left(\tilde{B}^{\left(\tau \right) } \right)}$ and ${{b_{n}}^{\left( k\right)}} {\left( q\right)}$ , 
   	\[\tilde{B}^{\left(\tau \right) }{{b_{n}}^{\left( k\right)}} {\left( q\right)}={e^{\left(k\right)}} \in c_{0} ,\] 
   	 Let  $x\in {\left[ c_{0} , \Delta^{\left(\tau \right)} , p \right]}_{e,r} $, then
  $ x^{\left[s \right]} =\sum _{k=0}^{s} \mu_{k}{\left( q\right)} {{b}^{\left( k\right)}}{\left( q\right)}$ for an integer $s \geq 0$.
  
   	By applying ${\tilde{B}^{\left(\tau \right) }}$ we get \ $ {\tilde{B}^{\left(\tau \right)}} x^{\left[s \right]} =\sum _{k=0}^{s} \mu _{k}{\left( q\right)} {\tilde{B}^{\left(\tau \right) }} {{b}^{\left( k\right)}}{\left( q\right)}$
  \[=\sum _{k=0}^{s} {\mu_{k}{\left( q\right)}}{e^{\left(k\right)}} = { \left(\left(  \tilde{B}^{\left({\tau}\right)}\right)x\right)}_k e^{\left(k\right)} \] and
 \begin{multline}
  	{\tilde{B}^{\left(\tau \right)} \left( {x-x^{\left[s\right]}}\right)}_{r} =\left\{\begin{array}{l}  {0, \qquad\qquad \qquad \qquad if\, 0\le r \le s } \\ {{\left(\left(  \tilde{B}^{\left({\tau}\right)}\right)x\right)}_k, \qquad \qquad  if\, r>s} \end{array}\right. ; for \ r,s \in {N_{0}}
  \end{multline} 
 \noindent For $\epsilon > 0$ there exist an integer $m_{o}$ s.t. 
   	\[{\mathop{\sup }\limits_{r \geq s}}\left|{ \left(\left(  \tilde{B}^{\left({\tau}\right)}\right)x\right)}_{r}\right| ^{\frac{p_{k}}{M}} < {\frac{\epsilon}{2}} \ for \ all \ s \geq {m_0} .\] 
   	Hence \[{g}_{\tilde{B}} \left( {x-x^{\left[s \right]}} \right) ={\mathop{\sup }\limits_{r \geq s}}\left|{ \left(\left(  \tilde{B}^{\left({\tau}\right)}\right)x\right)}_{r}\right| ^{\frac{p_{k}}{M}} < {\frac{\epsilon}{2}} < {\epsilon} \ , for \ all \ s\geq {m_0}. \]
   \noindent Assume that $x=\sum _{k}\eta _{k} {\left( q\right) } {b}^{\left(k\right)} {\left( q\right) }$. Since the linear mapping $F$ from $ {\left[ c_{0} , \Delta^{\left(\tau \right)} , p \right]}_{er} $ to $ {c_{0} \left( p\right)}$ is continuous we have,
 \[{ \left(\left(  \tilde{B}^{\left({\tau}\right)}\right)x\right)}_{k}=\sum _{k}\eta _{k} {\left( q\right)} \left( {\left(   \tilde{B}^{\left({\tau}\right)}\right)} {b}^{\left(k\right)} {\left( q\right)}\right)_n \]
 \[ =\sum _{k}\eta _{k} {\left( q\right)} e^{\left(k\right)} = \eta _{n} {\left( q\right)} \ . \]
   	This contradicts to our assumption that $ { \left(\left(  \tilde{B}^{\left({\tau}\right)}\right)x\right)}_{k} =\mu_{k} {\left( q\right)} $ for each $k\in {N_{0}}$. 
   	Thus, the representation is unique.
   	\noindent $(ii)$ The proof as it is similar to previous one.
   \end{proof}
  
    \section{${\alpha-}, \beta-$ and $\gamma-$ duals}
    
   Here we determine  $\alpha -,\beta -$ and $\gamma - $ duals  of  ${\left[ c_{0} , \Delta^{\left(\tau \right)} , p \right]}_{er}$ , ${\left[ c , \Delta^{\left(\tau \right)} , p \right]}_{er}$ and
   ${\left[ l_{\infty} , \Delta^{\left(\tau \right)} , p \right]}_{er}$.
  We define \[S(X,Y)=\left\{u=\left(u\right)_{k} \in \omega :ux=\left(u_{k} x_{k} \right)\in Y,whenever\, x=\left(x_{k} \right)\in X\right\}\] as the multiplier sequence space for any two sequence spaces  X and Y.
   Let $\alpha -,\beta -$ and $\gamma -$duals be denoted by 
   \noindent $X^{\alpha } =S\left(X,l_{1} \right),\, X^{\beta } =S\left(X,cs\right),\,  X^{\gamma } =S\left(X,bs\right)$ respectively.
   
   \noindent Throughout the collection of all finite subsets of $ \mathbb{N}$ is denoted by $\kappa $ . We consider $K\in \kappa$.
   \begin{l1}\cite{gros}\label{lm5}
   	\noindent Let $A=\left(a_{n,k} \right)$ be an infinite matrix. Then,
   	\begin{enumerate}
   	\item  $A\in \left(l_{\infty}{\left(p \right) }, l{\left(q \right) } \right)$ iff \begin{equation}
   			{\mathop{\sup }\limits_{k\in \kappa }} \sum _{n}{\left|\sum _{k\in K}a_{nk}{B^{\frac{1}{p_k}}}  \right|}^{q_n} <\infty ,\qquad for\, all\, integers\, B>1 and\, {q_n}\geq1 for\, all\, n ;
   		\end{equation}    
   	\item  $A\in  \left(l_{\infty}{\left(p \right) }, l_{\infty}{\left(q \right) } \right)$ iff  \begin{equation} \label{4.5} 
   			{\mathop{\sup } \limits_{n\in N}} \left( \sum _{k}\left|a_{nk} \right| {B^{\frac{1}{p_k}}} \right)^{q_n}  <\infty ,\qquad  ; 
   		\end{equation} 
   	\item  $A\in \left(l_{\infty}{\left(p \right) }, c{\left(q \right) } \right)$ and $q=\left({q_n} \right) $ be a bounded sequence of strictly positive real numbers iff \begin{equation} \label{4.6} 
   			{\mathop{\sup }\limits_{n\in N}} \sum _{k}\left|a_{nk} \right| {B^{\frac{1}{p_k}}}<\infty ,\qquad for\, all\, B>1 ,
   	\end{equation} exists $\left({\tau_{k}} \right)\subset R $  such that \   
   	${\mathop{\lim } \limits_{n\to \infty }}\left( \sum _{k}\left|  a_{nk}-{\tau_{k}} \right|{B^{\frac{1}{p_k}}} \right)^{q_n} =0 , for\, all\, B>1$;
    \item  $A\in \left(l_{\infty}{\left(p \right) }, c_{0}{\left(q \right) } \right)$ iff\[ {\mathop{\lim } \limits_{n\to \infty }}\left( \sum _{k}\left|  a_{nk} \right|{B^{\frac{1}{p_k}}} \right)^{q_n} =0 , \qquad for\, all\, B>1 , \] .
   \end{enumerate}
   \end{l1}
\begin{l1}\cite{gros}\label{lm6}
	\noindent Let $A=\left(a_{nk} \right)$ be an infinite matrix. Then,
	\begin{enumerate}
	\item  $A\in  \left(c_{0}{\left(p \right) }, l_{\infty}{\left(q \right) } \right)$ iff  \begin{equation} \label{4.7} 
			{\mathop{\sup } \limits_{n\in N}} \left( \sum _{k}\left|a_{nk} \right| {B^{\frac{-1}{p_k}}} \right)^{q_n}  <\infty ,\qquad for\, all\, B>1  ; 
		\end{equation} 
	\item  $A\in \left(c_{0}{\left(p \right) }, c{\left(q \right) } \right)$ iff \begin{equation} \label{4.8} 
		{\mathop{\sup }\limits_{n\in N}} \sum _{k}\left|a_{nk} \right| {B^{\frac{-1}{p_k}}}<\infty ,\qquad for\, all\, B>1 ,
	\end{equation}   
	\begin{equation} \label{4.9}  exists \left({\tau_{k}} \right)\subset R \,  such\, that\,
		{\mathop{\sup }\limits_{n\in N}}\sum _{k}\left|  a_{nk}-{\tau_{k}} \right|{M^{\frac{-1}{p_k}}}{B^{\frac{-1}{p_k}}} < \infty ,\\
		\end{equation} for all integers  $M,B>1$;
		\begin{equation} \label{5.1} exists \left({\tau_{k}} \right)\subset R \,   such\, that\,   
	{\mathop{\lim } \limits_{n\to \infty }} \sum _{k}\left|  a_{nk}-{\tau_{k}} \right|^{q_n} =0 , for\, all\, k \in N;
		\end{equation}
	\item  $A\in \left(c_{0}{\left(p \right) }, c_{0}{\left(q \right) } \right)$ iff
	\begin{equation} \label{5.2}  exists \left({\tau_{k}} \right)\subset R \,  such\, that\,
	{\mathop{\sup }\limits_{n\in N}}\sum _{k}\left|  a_{nk} \right|{M^{\frac{-1}{p_k}}}{B^{\frac{-1}{p_k}}} < \infty ,\\
\end{equation} for all integers  $M,B>1$;
\begin{equation} \label{5.3} exists \left({\tau_{k}} \right)\subset R \,   such\, that\,   
	{\mathop{\lim } \limits_{n\to \infty }} \sum _{k}\left|  a_{nk} \right|^{q_n} =0 , for\, all\, k \in N.
\end{equation}
	\end{enumerate}
\end{l1}
 \begin{l1}\cite{gros}\label{lm4}
 		\noindent Let $A=\left(a_{nk} \right)$ be an infinite matrix. Then,
 	\begin{enumerate}
 		\item  $A\in  \left(c{\left(p \right) }, l_{\infty}{\left(q \right) } \right)$ iff (16) holds and  \begin{equation} \label{4.2} 
 			{\mathop{\sup } \limits_{n\in N}}\left| \sum _{k}a_{nk} \right|^{q_n}  <\infty  ; 
 		\end{equation} 
 	\item  $A\in \left(c{\left(p \right) }, c{\left(q \right) } \right)$ iff (17), (18), (19) hold:
 \begin{equation} \label{5.4} exists \left({\tau} \right)\subset R \,   such\, that\,   
 	{\mathop{\lim } \limits_{n\to \infty }} \left| \sum _{k}  a_{nk}-{\tau} \right|^{q_n} =0 ;
 \end{equation}	
	\item  $A\in \left(c{\left(p \right) }, c_{0}{\left(q \right) } \right)$ iff (20), (21) hold and
\begin{equation} \label{5.5}
	{\mathop{\lim } \limits_{n\to \infty }} \left| \sum _{k}  a_{nk} \right|^{q_n} =0 .
\end{equation}	
\end{enumerate}		
 \end{l1}  
\begin{t1}
\noindent Let $0\ <{p_k} \leq H\ < \infty$   The ${\alpha-}, \beta-$ and $\gamma-$  duals of  ${\left[ c_{0} , \Delta^{\left(\tau \right)} , p \right]}_{er}$ , ${\left[ c , \Delta^{\left(\tau \right)} , p \right]}_{er}$ and
${\left[ l_{\infty} , \Delta^{\left(\tau \right)} , p \right]}_{er}$ are  the following sets 
\begin{tiny}
	\[D_{1}^{\left( \tau\right)}{\left( p\right) } =\\
	{\bigcap_{M>1}}\left\{a=\left(a_{k} \right)\in \omega :{\mathop{\sup }\limits_{k\in \tau }} \sum _{n}\left|\sum _{k\in K}\left[{\sum _{j=k}^{n-1}(-1)^{n-k} }  \left(\begin{array}{l} {j} \\ {k} \end{array}\right)\frac{\Gamma \left({-\tau} +1\right)}{\left(n-j\right)!\Gamma \left({-\tau} -n+j+1\right)} \ \frac{2^{k}  Q_{k}}{q_{j} }{a_k} \right] \right|M^{\frac{1}{p_k}} <\infty  \right\},\]
\[D_{2}^{\left( \tau\right)}{\left( p\right) } =
	{\bigcap_{M>1}}\left\{a=\left(a_{k} \right)\in \omega : \sum _{k}\left|\tilde{B}^{\left(\tau \right) }{\left(\frac{a_{k}}{q_{k}}\right)}{Q_k} \right|M^{\frac{1}{p_k}} <\infty   and \left(\frac{a_{k}}{q_{k}}{Q_k}M^{\frac{1}{p_k}} \right)\in {c_0}\right\} ,\] 
	\[D_{3}^{\left( \tau\right)}{\left( p\right) } =
	{\bigcap_{M>1}}\left\{a=\left(a_{k} \right)\in \omega :\sum _{k}\left|\tilde{B}^{\left(\tau \right)}{\left(\frac{a_{k}}{q_{k}}\right)}{Q_k} \right|M^{\frac{1}{p_k}} <\infty   and \left\lbrace\tilde{B}^{\left(\tau \right)}{\left(\frac{a_{k}}{q_{k}}\right)}{Q_k}  \right\rbrace \in {l_ {\infty}} \right\},\] 
	\[D_{4}^{\left( \tau\right)}{\left( p\right) } =\\
{\bigcup_{M>1}}\left\{a=\left(a_{k} \right)\in \omega :{\mathop{\sup }\limits_{k\in \tau }} \sum _{n}\left|\sum _{k\in K}\left[{\sum _{j=k}^{n-1}(-1)^{n-k} }  \left(\begin{array}{l} {j} \\ {k} \end{array}\right)\frac{\Gamma \left({-\tau} +1\right)}{\left(n-j\right)!\Gamma \left({-\tau} -n+j+1\right)} \ \frac{2^{k}  Q_{k}}{q_{j} }{a_k} \right] \right|M^{\frac{-1}{p_k}} <\infty  \right\},\]
\[D_{5}^{\left( \tau\right)}{\left( p\right) } =\\
{\bigcup_{M>1}}\left\{a=\left(a_{k} \right)\in \omega : \sum _{n}\left|\sum _{k}\left[{\sum _{j=k}^{n-1}(-1)^{n-k} }  \left(\begin{array}{l} {j} \\ {k} \end{array}\right)\frac{\Gamma \left({-\tau} +1\right)}{\left(n-j\right)!\Gamma \left({-\tau} -n+j+1\right)} \ \frac{2^{k}  Q_{k}}{q_{j} }{a_k} \right] \right| <\infty  \right\},\]	
\[D_{6}^{\left( \tau\right)}{\left( p\right) } =
{\bigcap_{M>1}}\left\{a=\left(a_{k} \right)\in \omega :\sum _{k}\left|\tilde{B}^{\left(\tau \right)}{\left(\frac{a_{k}}{q_{k}}\right)}{Q_k} \right|M^{\frac{-1}{p_k}} <\infty  \right\},\] 
\end{tiny}
where\begin{equation}
	\tilde{B}^{\left(\tau \right)}{\left(\frac{a_{k}}{q_{k}}\right)} ={\frac{2^{k} a_{k}}{q_{k}}+{\sum _{i=k+1}^{n}(-1)^{n-k}{a_i} } {\sum _{j=k}^{i}} \left(\begin{array}{l} {i} \\ {j} \end{array}\right)\frac{\Gamma \left({-\tau} +1\right)}{\left(j-i\right)!\Gamma \left({-\tau} -j+i+1\right)} \ \frac{2^{i}}{q_{j} }}
\end{equation}
\end{t1}  
\begin{proof}
	\noindent We prove it for ${\left[ l_{\infty} , \Delta^{\left(\tau \right)} , p \right]}_{er}$ . Considering $x=\left(x_{k} \right)$ as in \ref{3.2}, let $a=\left(a_{k} \right)\in \omega $ define
	\[\begin{array}{l} {a_{n} x_{n} =\sum _{i=0}^{k} \sum _{j=i}^{k}\left(-1\right)^{k-i}  \left(\begin{array}{l} {j} \\ {i} \end{array}\right)\frac{\Gamma \left(-{\tau} +1\right)}{\left(k-j\right)!\Gamma \left(-{\tau} -k+j+1\right)}  \frac{2^i Q_{i}}{q_{j}}{a_n} y_{i} } \\ {\, \, \, \, \, \, \, \, \, \, =\left(Uy \right)_{n} ,\qquad {\rm for}\, n\in N} \end{array},\]
	  where matrix $U =\left(u_{nk} \right)$ is defined as
	  \begin{multline}
	 	 {u_{nk}}=\left\{\begin{array}{l} {\sum _{j=k}^{n}(-1)^{n-k} }  \left(\begin{array}{l} {j} \\ {k} \end{array}\right)\frac{\Gamma \left({-\tau} +1\right)}{\left(n-j\right)!\Gamma \left({-\tau} -n+j+1\right)} \ \frac{2^{k}  Q_{k}}{q_{j}{a_n} } , { \qquad if\, 0\le k \le n } \\\ {\frac{2^{n} Q_{n}}{q_{n}}, \qquad \qquad \qquad \qquad \qquad \qquad\qquad\qquad \qquad if\, k=n} \\\ \   {0,\qquad \qquad \qquad \qquad \qquad \qquad\qquad \qquad \qquad \qquad if\, k>n} \end{array}\right. 
	 \end{multline}  
Therefore we conclude that $ax=\left(d_{n} x_{n} \right)\in l_{1} $ when $x=\left(x_{k} \right) \in {\left[ l_{\infty} , \Delta^{\left(\tau \right)} , p \right]}_{er}$	iff $Uy \in l_{1}$ as $y=\left({y_k} \right) \in {l_{\infty}}{\left( p\right)} $ .\\
By lemma $\left( \ref{lm5}\right) $ we conclude that $\left\lbrace {\left[ l_{\infty} , \Delta^{\left(\tau \right)} , p \right]}_{er} \right\rbrace ^{\alpha}=D_{1}^{\left( \tau\right)}{\left( p\right) } $.
\noindent Now,
\[\sum _{i=0}^{n}{a_{k} x_{k} =\sum _{k=0}^{n} {a_k}\left[ \sum _{i=0}^{k} \sum _{j=i}^{k}\left(-1\right)^{k-i}  \left(\begin{array}{l} {j} \\ {i} \end{array}\right)\frac{\Gamma \left(-{\tau} +1\right)}{\left(k-j\right)!\Gamma \left(-{\tau} -k+j+1\right)}  \frac{2^i Q_{i}}{q_{j}} y_{i} \right] }\]
\[ =\sum _{k=0}^{n} {y_k}{Q_k} \left[ {\frac{2^{k} a_{k}}{q_{k}}+{\sum _{i=k+1}^{n}(-1)^{n-k}{a_i} } {\sum _{j=k}^{i}} \left(\begin{array}{l} {i} \\ {j} \end{array}\right)\frac{\Gamma \left({-\tau} +1\right)}{\left(j-i\right)!\Gamma \left({-\tau} -j+i+1\right)} \ \frac{2^{i}}{q_{j} }} \right] \]
\[ =\sum _{k=0}^{n} {y_k}{Q_k} \tilde{B}^{\left(\tau \right)}{\left(\frac{a_{k}}{q_{k}}\right)} = \left(Vy \right)_{n},\]
  where matrix $V =\left(v_{nk} \right)$ is defined as
      \begin{multline}
     	{v_{nk}}=\left\{\begin{array}{l} {\tilde{B}^{\left(\tau \right)}{\left(\frac{a_{k}}{q_{k}}\right)}{Q_k}}, {\qquad \qquad \qquad \qquad if\, 0\le k \le n} \\{\frac{2^{n} Q_{n}}{q_{n}}, \qquad \qquad \qquad \qquad  \qquad \qquad  if\, k=n} \\ {0,\qquad \qquad \qquad \qquad \qquad \qquad \qquad if\, k>n} \end{array}\right. 
        \end{multline}  
Therefore we deduce that $ax=\left(d_{n} x_{n} \right)\in cs $ when $x=\left(x_{k} \right) \in {\left[ l_{\infty} , \Delta^{\left(\tau \right)} , p \right]}_{er}$	iff $Vy \in c$ as $y=\left({y_k} \right) \in {l_{\infty}}{\left( p\right)} $ .\\
By using lemma $\left( \ref{lm5}\right) $ with $ q=q_{n}=1 $ we conclude that $\left\lbrace {\left[ l_{\infty} , \Delta^{\left(\tau \right)} , p \right]}_{er} \right\rbrace ^{\beta}=D_{2}^{\left( \tau\right)}{\left( p\right) } $.
Similarly by using  lemma $\left( \ref{lm5}\right) $ with $ q=q_{n}=1 $ for all n, we conclude that $\left\lbrace {\left[ l_{\infty} , \Delta^{\left(\tau \right)} , p \right]}_{er} \right\rbrace ^{\gamma}=D_{3}^{\left( \tau\right)}{\left( p\right) } $.Hence the theorem proved and the duals of other spaces can be obtained in a similar manner using lemma $\left( \ref{lm6}\right) $ and lemma $\left( \ref{lm4}\right) $.
\end{proof}

     \end{document}